\documentclass{amsart}
\usepackage{amsfonts}
\usepackage{graphicx}

\newtheorem{theorem}{Theorem}[section]
\newtheorem*{theorem A}{Theorem A}
\newtheorem*{theorem B}{N\"olker's Theorem}

\newtheorem{corollary}{Corollary}[section]

\theoremstyle{remark}
\newtheorem{remark}{Remark}[section]
\theoremstyle{remark}

\theoremstyle{definition}

\newtheorem{example}{Example}[section]
\numberwithin{equation}{section}
\def\({\left ( }
\def\){\right )}
\def\<{\left < }
\def\>{\right >}

\input{tcilatex}

\begin{document}
\title{ Certain classifications on surfaces of revolution in a
semi-isotropic space}
\author{Muhittin Evren Aydin}
\address{Department of Mathematics, Faculty of Science, Firat University,
Elazig, 23119, Turkey}
\email{meaydin@firat.edu.tr}
\thanks{}
\subjclass[2000]{53A35, 53B25, 53B30, 53C42.}
\keywords{Semi-isotropic space, Surface of revolution, Laplace
operator.}

\begin{abstract}
A semi-isotropic space is a real affine 3-space endowed
with the non-degenerate metric $dx^{2}-dy^{2}.$ The main purpose of this paper is to describe
the surfaces of revolution in the semi-isotropic space that satisfy some
equations in terms of the position vector and the Laplace operators with
respect to the first and the second fundamental forms.
\end{abstract}

\maketitle

\section{introduction}

The  \textit{isotropic 3-space} $\mathbb{I}^{3}$ is one of the Cayley-Klein 3-spaces
and can be defined in the real projective 3-space. One can also appear as the affine 3-space $\mathbb{R}^{3}$  equipped with the semi-norm (see \cite{1,2,28,32,33})%
\begin{equation*}
\left\Vert u\right\Vert _{\mathbb{I}}=\sqrt{\left( u_{1}\right) ^{2}+\left(
u_{2}\right) ^{2}},\text{ }u=\left( u_{1},u_{2},u_{3}\right) \in \mathbb{I}%
^{3}.
\end{equation*}

The isotropic geometry has remarkable applications in Image Processing, architectural design and microeconomics, see \cite%
{10,13,20}, \cite{29}-\cite{31}.

The fundamentals of curves and surfaces in $\mathbb{I}^{3}$ can be found in
H. Sachs' monograph \cite{32}. For further studies of these in $\mathbb{I}^{3}$
we refer to \cite{6,14,18,23,24,27}.

A \textit{semi-isotropic 3-space} $\mathbb{SI}^{3}$ is the
product of the Lorentz-Minkowsi 2-space $%
\mathbb{E}
_{1}^{2}$ and the isotropic line equipped with a degenerate parabolic
distance metric. More precisely, $\mathbb{SI}^{3}$ is the affine 3-space $\mathbb{R}^{3}$ 
endowed with the semi-norm%
\begin{equation*}
\left\Vert u\right\Vert =\sqrt{\left\vert \left( u_{1}\right) ^{2}-\left(
u_{2}\right) ^{2}\right\vert },\text{ }u=\left( u_{1},u_{2},u_{3}\right) \in 
\mathbb{SI}^{3}.
\end{equation*}

The local theory of non-null curves and surfaces in $\mathbb{SI}^{3}$ was
recently stated in \cite{3}. This has caused that many (open or solved) problems
of classic differential geometry can be treated to $\mathbb{SI}^{3}$. For
example, the surfaces of revolution in $\mathbb{SI}^{3}$ with constant
curvature were classified in \cite{3}.

In this paper we aim to present the surfaces of revolution in $%
\mathbb{SI}^{3}$ that satisfy certain conditions in terms of the coordinate
functions of the position vector and the Laplace operator.

These are natural, being related to the so-called \textit{%
submanifolds of finite type}, introduced by B.-Y. Chen in the late 1970's
(see \cite{7}-\cite{9},\cite{11}).

Let $r$ be an isometric immersion of a Riemannian manifold $M$ into the
Euclidean $n-$space $\mathbb{E}^{n}$ and $\bigtriangleup $ denote the
Laplace operator of $M$. Then it is said to be of \textit{finite type }if
its position vector field can be expressed as%
\begin{equation*}
r=c+r_{0}+r_{1}+...+r_{k},
\end{equation*}%
where $c\in \mathbb{R}^{n}$ and $x_{i}$ non-constant $\mathbb{R}^{n}-$valued
maps such that%
\begin{equation*}
\bigtriangleup r_{i}=\lambda _{i}r_{i},\text{ }1\leq i\leq k.
\end{equation*}%
If $\lambda _{1},...,\lambda _{k}$ are mutually different, then the
immersion is said to be of $k-$\textit{type. }If one of $\lambda
_{1},...,\lambda _{k}$ is zero, then the immersion is said to be of \textit{%
null }$k-$\textit{type.}

Additionaly,  O. Garay classified the
surfaces of revolution in $\mathbb{R}^{3}$
satisfying the following relation (\cite{15,16})%
\begin{equation}
\bigtriangleup r_{i}=\lambda _{i}r_{i},  \tag{1.1}
\end{equation}%
where $r^{i}$ are component functions of the position
vector. This idea was generalized to $\mathbb{R}_{1}^{3}$ in \cite{5,12,17}.

The results of this paper inherently resemble to these of existing ones in
literature; however, the discrepancies also arise. For instance, in Theorem 4.2, we prove
that a surface of revolution in $\mathbb{SI}^{3}$ that satisfy $\left(
1.1\right) $ is formed by rotating the graph of a Bessel function.

Such functions which correspond to the solutions of Bessel's equation are
named due to the German mathematician and astronomer Frederic Wilhelm Bessel
(1784--1846) who first used them to analyze planetary orbits. The Bessel functions are widely used in classical physics, elasticity theory, heat conduction theory etc. (see \cite{21} ). We give a
brief review of the Bessel functions in \S 3.

Morever we state that the surfaces of revolution in $\mathbb{SI}^{3}$ that
satisfy%
\begin{equation*}
\bigtriangleup ^{II}r_{i}=\lambda_{i}r_{i}
\end{equation*}%
are only (s-i)-minimal ones, where $\bigtriangleup ^{II}$ is the Laplace
operator with respect to the second fundamantal form (see Theorem 4.3, Corollary4.1).

\section{Preliminaries}

We provide the basics of semi-isotropic space from \cite%
{3}.

The \textit{semi-isotropic geometry} is based on the six-parameter group of affine
transformations%
\begin{equation}
\left\{ 
\begin{array}{l}
\mathbf{x}^{\prime }=a_{1}+\mathbf{x}\cosh a_{2}+\mathbf{y}\sinh a_{2}, \\ 
\mathbf{y}^{\prime }=a_{3}+\mathbf{x}\sinh a_{2}+\mathbf{y}\cosh a_{2}, \\ 
\mathbf{z}^{\prime }=a_{4}+a_{5}\mathbf{x}+a_{6}\mathbf{y}+\mathbf{z,}\text{ 
}a_{i}\in 
\mathbb{R}
,\text{ }i=1,...,6.%
\end{array}%
\right.  \tag{2.1}
\end{equation}%
These are called \textit{semi-isotropic congruence transformations} or%
\textit{\ }(\textit{s-i)-motions. }

We notice that the (s-i)-motions are a composition of a Lorentzian motion in 
$\mathbf{xy}-$plane and an affine shear transformation in $\mathbf{z}$%
-direction. For detailed properties of the Lorentzian motions, see \cite%
{22,26}.

The\textit{\ semi-isotropic scalar product }between two vectors $u=\left(
u_{1},u_{2},u_{3}\right) $ and $v=\left( v_{1},v_{2},v_{3}\right) \in 
\mathbb{SI}^{3}$\textit{\ }is given by%
\begin{equation}
\left\langle u,v\right\rangle =\left\{ 
\begin{array}{ll}
u_{3}v_{3}, & \text{if }u_{1}=u_{2}=v_{1}=v_{2}=0\text{,} \\ 
u_{1}v_{1}-u_{2}v_{2}, & \text{otherwise.}%
\end{array}%
\right.  \notag
\end{equation}

The \textit{vector product }in the sense of semi-isotropic space is%
\begin{equation*}
u\times v=%
\begin{vmatrix}
\left( 1,0,0\right) & \left( 0,-1,0\right) & \left( 0,0,0\right) \\ 
u_{1} & u_{2} & u_{3} \\ 
v_{1} & v_{2} & v_{3}%
\end{vmatrix}%
.
\end{equation*}%
For some vector $w=\left( w_{1},w_{2},w_{3}\right) \in \mathbb{SI}^{3},$
denote $\tilde{w}$ its canonical projection onto $%
\mathbb{E}
_{1}^{2},$ i.e. $\tilde{w}=\left( w_{1},w_{2},0\right) .$ It is then seen that 
\begin{equation*}
\left\langle u\times v,w\right\rangle =\det \left( u,v,\tilde{w}\right) .
\end{equation*}

We call the vector $u\neq 0$ \textit{isotropic }if $\tilde{u}$ is a zero vector.
If $u=0$ or $\tilde{u}\neq 0$ it is called \textit{non-isotropic.}

A non-isotropic vector $u\in \mathbb{SI}^{3}$ is respectively called \textit{%
spacelike}, \textit{timelike} and \textit{null }(or\textit{\ lightlike}) if $%
\left\langle u,u\right\rangle >0$ or $u=0,$ $\left\langle u,u\right\rangle
<0 $ and $\left\langle u,u\right\rangle =0$ $\left( u\neq 0\right) .$

The \textit{null-cone }and the \textit{timelike-cone} of $\mathbb{SI}^{3}$
are respectively given by%
\begin{equation*}
\mathcal{C}=\left\{ \left. \left( x,y,z\right) \in \mathbb{SI}%
^{3}\right\vert x^{2}-y^{2}=0\right\} -\left\{ 0\in \mathbb{SI}^{3}\right\} .
\end{equation*}%
and%
\begin{equation*}
\mathcal{T}=\left\{ \left. \left( x,y,z\right) \in \mathbb{SI}%
^{3}\right\vert x^{2}-y^{2}<0\right\} .
\end{equation*}

The \textit{semi-isotropic} \textit{angle} between two timelike vectors $%
u,v\in \mathbb{SI}^{3}$ is defined as 
\begin{equation*}
\left\langle u,v\right\rangle =-\left\Vert u\right\Vert \left\Vert
v\right\Vert \cosh \phi .
\end{equation*}

Note that all isotropic vectors are orthonogal to non-isotropic ones. Also,
two non-isotropic vectors $u,v$ in $\mathbb{SI}^{3}$ are orthonogal if $%
\left\langle u,v\right\rangle =0.$

Let $\alpha \left( s\right) $ be a regular curve in $\mathbb{SI}^{3},$ i.e. $%
\alpha ^{\prime }\left( s\right) =\frac{d\alpha }{ds}\neq 0$ for all $s.$
Then it is said to be \textit{admissible} if $\alpha \left( s\right) $ has
no isotropic tangent vector, i.e. $\tilde{\alpha}^{\prime }\left( s\right)
\neq 0.$ An admissible curve $\alpha \left( s\right) $ in $\mathbb{SI}^{3}$
is called \textit{spacelike} (resp. \textit{timelike}, \textit{null})\textit{%
\ }if $\alpha ^{\prime }\left( s\right) $ is spacelike (resp. timelike,
null) for all $s.$

\subsection{Planes, circles and spheres in $\mathbb{SI}^{3}$}

The lines in $\mathbf{z}-$direction are called \textit{isotropic lines}. The
planes containing an isotropic line are called \textit{isotropic planes}.
Other planes are \textit{non-isotropic}.

In the non-isotropic planes the Lorentzian metric is basically used. Let $%
\Gamma $ be a non-isotropic plane and $c_{0}=\left( x_{0},y_{0},z_{0}\right) 
$ a fixed point in $\Gamma $. The iso-distance set of $c_{0}$ in $\Gamma $
is 
\begin{equation*}
\left\{ p=\left( x,y,z\right) \in \mathbb{SI}^{3}:\left\Vert
p-c_{0}\right\Vert =r,r\in 
\mathbb{R}
^{+}\right\} .
\end{equation*}%
The projection of such a set onto $\mathbf{xy}-$plane is a rectangular
hyperbola. We call it a \textit{semi-isotropic circle} (or \textit{%
(s-i)-circle)} \textit{of hyperbolic type}.

In the isotropic planes we have an isotropic metric. The\textit{\
(s-i)-circle} \textit{of parabolic type} is a parabola with $\mathbf{z}-$%
axis lying in an isotropic plane. An (s-i)-circle of parabolic type is not
the iso-distance set of a fixed point.

We have two types of semi-isotropic spheres. One is the \textit{%
semi-isotropic sphere }(or\textit{\ (s-i)-sphere})\textit{\ of cylindrical
type} which is the set of all points $p\in \mathbb{SI}^{3}$ with $\left\Vert
p-c_{0}\right\Vert =r.$ This sphere is a right Lorentz hyperbolic cylinder
with $\mathbf{z}-$parallel rulings from the Lorentz-Minkowski perspective.

Other one is the \textit{(s-i)-sphere of parabolic type,}%
\begin{equation*}
z=r\left( x^{2}-y^{2}\right) +ax+by+c,\text{ }r>0.
\end{equation*}%
\textit{\ }These are indeed hyperbolic paraboloids.

The intersections of these (s-i)-spheres with non-isotropic (resp.
isotropic) planes are (s-i)-circles of hyperbolic type (resp. of parabolic
type).

\subsection{Spacelike and timelike surfaces in $\mathbb{SI}^{3}$}

Let $M$ be an admissible surface immersed in $\mathbb{SI}^{3},$ i.e.,
without isotropic tangent planes. Denote $g$ the metric on $M$ induced
from $\mathbb{SI}^{3}.$ The surface $M$ is said to be \textit{spacelike}
(resp. \textit{timelike}, \textit{null})\textit{\ }if $g$ is positive
definite (resp. a metric with index 1, degenerate).

The following contains the formulas for only spacelike and timelike
admissible surfaces in $\mathbb{SI}^{3}.$

Assume that $M$ has a local parameterization of the form%
\begin{equation*}
r:D\subseteq \mathbb{R}^{2}\longrightarrow \mathbb{SI}^{3}:\text{ }\left(
u,v\right) \longmapsto \left( x\left( u,v\right) ,y\left( u,v\right)
,z\left( u,v\right) \right) .
\end{equation*}

The components of the first fundamental form $I$ with respect to the basis $%
\left\{ r_{u},r_{v}\right\} $ are 
\begin{equation}
E=\left\langle r_{u},r_{u}\right\rangle ,\text{ }F=\left\langle
r_{u},r_{v}\right\rangle ,\text{ }G=\left\langle r_{v},r_{v}\right\rangle , 
\notag
\end{equation}%
where $r_{u}=\partial r/\partial u.$ $M$ is timelike (spacelike) when $%
W=EG-F^{2}<0$ $(>0).$

Denote $\bigtriangleup$ the Laplace operator of $M$. For a smooth function $\psi :M\longrightarrow 
\mathbb{R}
,$ $\left( u,v\right) \longmapsto \psi \left( u,v\right) $ its Laplacian is
defined as 
\begin{equation}
\bigtriangleup \psi =-\frac{1}{\sqrt{\left\vert W\right\vert }}\left\{ 
\frac{\partial }{\partial u}\left( \frac{G\psi _{u}-F\psi _{v}}{\sqrt{%
\left\vert W\right\vert }}\right) -\frac{\partial }{\partial v}\left( \frac{%
F\psi _{u}-E\psi _{v}}{\sqrt{\left\vert W\right\vert }}\right) \right\} . 
\tag{2.2}
\end{equation}

The unit normal vector field of $M$ is the isotropic vector $\left(
0,0,1\right) $ since it is perpendicular to all non-isotropic vectors.

The components of the second fundamental form $II$ are given by%
\begin{equation}
L=\frac{\det \left( r_{uu},r_{u},r_{v}\right) }{\sqrt{\left\vert
W\right\vert }},\text{ }M=\frac{\det \left( r_{uv},r_{u},r_{v}\right) }{%
\sqrt{\left\vert W\right\vert }},\text{ }M=\frac{\det \left(
r_{vv},r_{u},r_{v}\right) }{\sqrt{\left\vert W\right\vert }}  \notag
\end{equation}%
for $r_{uu}=\frac{\partial ^{2}r}{\partial u^{2}},$ etc.

Thus the \textit{semi-relative curvature} and the \textit{semi}-\textit{%
isotropic mean curvature} of $M$ are respectively defined by%
\begin{equation}
K=-\varepsilon \frac{LN-M^{2}}{W}\text{ and }H=-\varepsilon \frac{EN-2FM+GL}{%
2W},  \notag
\end{equation}%
where $\varepsilon =sgn\left( W\right) $.

Assume that $M$ has no parabolic points, i.e. $w=LN-M^{2}\neq 0.$ Then the
Laplace operator with respect to $II$ is given by%
\begin{equation}
\bigtriangleup ^{II}\psi =-\frac{1}{\sqrt{\left\vert w\right\vert }}\left\{ 
\frac{\partial }{\partial u}\left( \frac{N\psi _{u}-M\psi _{v}}{\sqrt{%
\left\vert w\right\vert }}\right) -\frac{\partial }{\partial v}\left( \frac{%
M\psi _{u}-L\psi _{v}}{\sqrt{\left\vert w\right\vert }}\right) \right\} . 
\tag{2.3}
\end{equation}

\section{Bessel's differential equation}

This section is devoted to recall the solutions of the following linear
second-order ordinary differential equation (ODE)%
\begin{equation}
x^{2}\frac{d^{2}y}{dx^{2}}+x\frac{dy}{dx}+\left( x^{2}-p^{2}\right) y=0 
\tag{3.1}
\end{equation}%
for a real positive constant $p$ (see \cite{4,19,25}). This is called 
\textit{Bessel's equation of order} $p$. If the sign of the term $\left(
x^{2}-p^{2}\right) $ is negative, then it is called \textit{Bessel's
modified equation of order }$p.$ Since the ODE $\left( 3.1\right) $ has a
singular point at $x=0,$ a series solution for $\left( 3.1\right) $ via the
method of Frobenius can be construct. Before this, it is useful to express
the gamma function and the Pockhammer symbol.

The \textit{gamma function} is defined as%
\begin{equation*}
\Gamma \left( x\right) =\int_{0}^{\infty }e^{-q}q^{x-1}dq,\text{ }x>0.
\end{equation*}%
If $x=n$ is a positive integer, then it can be seen that%
\begin{eqnarray*}
\Gamma \left( n+1\right) &=&n\Gamma \left( n\right) =n\left( n-1\right)
\Gamma \left( n-1\right) =...= \\
&=&n\left( n-1\right) ...2.1=n!.
\end{eqnarray*}

The \textit{Pockhammer symbol }provides a simplicity for the product of too
many terms and is given by%
\begin{equation*}
\left( k\right) _{n}=k\left( k+1\right) \left( k+2\right) ...\left(
k+n-1\right) .
\end{equation*}%
In general $\left( k\right) _{0}=1$.

The indicial equation for $\left( 3.1\right) $ implies the roots $%
r_{1,2}=\pm p.$ When $2p$ is not an integer, the method of Frobenius gives
two linearly independent solutions as follows%
\begin{equation*}
y_{1,2}\left( x\right) =\frac{x^{\pm p}}{2^{\pm p}\Gamma \left( 1\pm
p\right) }\sum_{n=0}^{\infty }\frac{\left( -1\right) ^{n}x^{2n}}{%
2^{2n}\left( 1\pm p\right) _{n}n!}.
\end{equation*}%
These are called \textit{Bessel functions of order }$\pm p$ denoted by $%
J_{^{\pm p}}.$ Thus the general solution of $\left( 3.1\right) $ becomes%
\begin{equation*}
y\left( x\right) =c_{1}J_{p}\left( x\right) +c_{2}J_{-p}\left( x\right) ,%
\text{ }c_{1},c_{2}\in 
\mathbb{R}
.
\end{equation*}

In case $p=0,$ by applying again the method of Frobenius, we obtain that the
first solution is%
\begin{equation*}
J_{0}\left( x\right) =\sum_{n=0}^{\infty }\frac{\left( -1\right) ^{n}x^{2n}}{%
2^{2n}\left( n!\right) ^{2}}
\end{equation*}%
and second%
\begin{equation*}
Y_{0}\left( x\right) =\frac{2}{\pi} \left\{ \left( \ln \frac{x}{2}+\gamma \right)J_{0}-\sum_{n=0}^{\infty }\frac{%
\left( -1\right) ^{n}\phi \left( n\right) x^{2n}}{2^{2n}\left( n!\right) ^{2}%
} \right\},
\end{equation*}%
where $\phi \left( n\right) =\sum_{m=1}^{n}\frac{1}{m}$ and $\gamma$  Euler-Mascheroni constant.  Such functions are
respectively called the \textit{Bessel function of the first kind of order
zero }and \textit{Weber's Bessel function of order zero}, see Fig. 1(a).

For the Bessel's modified equation of order zero, the general solution is%
\begin{equation*}
y\left( x\right) =c_{1}I_{0}\left( x\right) +c_{2}K_{0}\left( x\right) ,%
\text{ }c_{1},c_{2}\in 
\mathbb{R}
,
\end{equation*}%
where%
\begin{equation*}
I_{0}\left( x\right) =\sum_{n=0}^{\infty }\frac{x^{2n}}{2^{2n}\left(
n!\right) ^{2}}
\end{equation*}%
and%
\begin{equation*}
K_{0}\left( x\right) =-\left( \ln \frac{x}{2}+\gamma \right) I_{0}\left(
x\right) +\sum_{n=0}^{\infty }\frac{\phi \left( n\right) x^{2n}}{%
2^{2n}\left( n!\right) ^{2}},
\end{equation*}%
which are respectively called \textit{modified Bessel functions of first and
second kind of order }$0$, see Fig 1(b).%

\bigskip

\begin{minipage}{0.5\textwidth}
\includegraphics[scale=0.48]{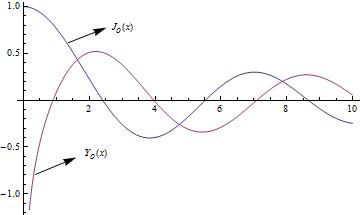} 
\verb Figure1(a).  Bessel functions of order zero, $x\in[0,10]$. 
\end{minipage}
\begin{minipage}{0.45\textwidth}
\includegraphics[scale=0.45]{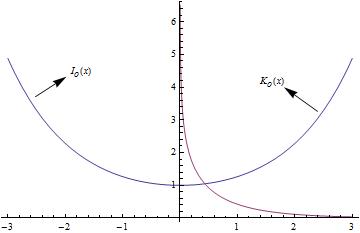}
\verb Figure1(b).  Modified Bessel functions of order zero, $x\in[-3,3]$.
\end{minipage}

\section{Classifications on surfaces of revolution in $\mathbb{SI}^{3}$}

Let $\gamma \left( u\right) =\left( 0,u,f\left( u\right) \right) $ (resp. $%
\gamma \left( u\right) =\left( u,0,f\left( u\right) \right) $ be a timelike
(resp. spacelike) admissible curve in $\mathbb{SI}^{3}$. By rotating them
around $\mathbf{z}-$axis via the transformations $\left( 2.1\right) $ we
derive%
\begin{equation}
r\left( u,v\right) =\left( u\sinh v,u\cosh v,f\left( u\right) \right) 
\tag{4.1}
\end{equation}%
and 
\begin{equation}
r\left( u,v\right) =\left( u\cosh v,u\sinh v,f\left( u\right) \right) . 
\tag{4.2}
\end{equation}%
These are called \textit{surfaces of revolution} in $\mathbb{SI}^{3}.$ Also
we call $\gamma $ \textit{profile curve. } Remark that such surfaces of $\mathbb{SI}^{3}$ are timelike since $EG-F^{2}=-u^{2}.$

The (s-r)-curvature $K$ and (s-i)-mean curvature $H$ of these surfaces in $%
\mathbb{SI}^{3}$ are%
\begin{equation*}
K=\frac{f^{\prime }f^{\prime \prime }}{u}\text{ and }H=\frac{1}{2}\left( 
\frac{f^{\prime }}{u}+f^{\prime \prime }\right) ,
\end{equation*}%
where $f^{\prime }\left( u\right) =\frac{df}{du}$ and so on.

The following classifies the surfaces of revolution in $\mathbb{SI}^{3}$
with $K=const.$ and $H=const.$

\begin{theorem} \cite{3} Let $M$ be a surface of
revolution in $\mathbb{SI}^{3}.$ Then the following statements hold:%

(i) $M$ has nonzero constant (s-r)-curvature $K_{0}$ 
 if and only if its profile curve is of the form%
\begin{equation*}
\left\{ 
\begin{array}{l}
\left( u,0,f\left( u\right) \right) \text{ or }\left( 0,u,f\left( u\right)
\right) \\ 
f\left( u\right) =\frac{u}{2}\psi \left( u\right) +\frac{c_{1}}{2\sqrt{K_{0}}%
}\ln \left\vert 2\sqrt{K_{0}}\left( \sqrt{K_{0}}x+\psi \left( u\right)
\right) \right\vert , \\ 
\psi \left( u\right) =\sqrt{c_{1}+K_{0}u^{2}},c_{1},c_{2}\in 
\mathbb{R}
.%
\end{array}%
\right.
\end{equation*}%
In the particular case $K=0,$ $f$ becomes a
linear function.

(ii)$M$ has constant (s-im)-curvature $H_{0}$ if
and only if its profile curve is given by%
\begin{equation*}
\left\{ 
\begin{array}{l}
\left( u,0,f\left( u\right) \right) \text{ or }\left( 0,u,f\left( u\right)
\right) \\ 
f\left( u\right) =\frac{H_{0}}{2}u^{2}+c_{1}\ln u+c_{2},\text{ }%
c_{1},c_{2}\in 
\mathbb{R}
.%
\end{array}%
\right.
\end{equation*}
\end{theorem}

Now let $M$ be a surface of revolution in $\mathbb{SI}^{3}$ given by $\left(
4.1\right) $. Then we have%
\begin{equation}
r_{1}\left( u,v\right) =u\sinh v,\text{ }r_{2}\left( u,v\right) =u\cosh v,%
\text{ }r_{3}\left( u,v\right) =f\left( u\right) .  \tag{4.3}
\end{equation}%
Denote $\bigtriangleup$ the Laplacian of $M$. It
follows from $\left( 2.2\right) $ and $\left( 4.3\right) $ that%
\begin{equation}
\bigtriangleup r_{1} =\bigtriangleup r_{2} =0  \tag{4.4}
\end{equation}%
and%
\begin{equation}
\bigtriangleup  r_{3} =-f^{\prime \prime
}-\frac{1}{u}f^{\prime }.  \tag{4.5}
\end{equation}%
Assume that $M$ is not (s-i) minimal and $\bigtriangleup r_{i}=\lambda
_{i}r_{i},$ $i=1,2,3.$ Then from $\left( 4.4\right) $ and $\left( 4.5\right) 
$ we find that $\lambda _{1}=\lambda _{2}=0$ and%
\begin{equation}
f^{\prime \prime }+\frac{1}{u}f^{\prime }+\lambda _{3}f=0, \text{ } \lambda _{3} \neq 0. \tag{4.6}
\end{equation}%
If $\lambda _{3}>0$ (resp. $<0$) in $\left( 4.6\right) ,$ then it becomes a
Bessel's (resp. modified) ODE of order zero. The general solutions for $%
\left( 4.6\right) $ turn to%
\begin{equation*}
f\left( u\right) =c_{1}J_{0}\left( \sqrt{\lambda _{3}}u\right) +c_{2}Y_{0}\left( 
\sqrt{\lambda _{3}}u\right)
\end{equation*}%
and (in case $\lambda _{3}<0$)%
\begin{equation*}
f\left( u\right) =c_{1}I_{0}\left( \sqrt{-\lambda _{3}}u\right) +c_{2}K_{0}\left( 
\sqrt{-\lambda _{3}}u\right) ,\text{ }c_{1},c_{2}\in 
\mathbb{R}
.
\end{equation*}

\begin{theorem} Let $M$\ be a surface of revolution in $\mathbb{SI}^{3}$%
\ given by $\left( 4.1\right) .$\ If $\bigtriangleup r_{i}=\lambda
_{i}r_{i}$\ then $\left( \lambda _{1},\lambda _{2},\lambda _{3}\right)
=\left( 0,0,\lambda \neq 0\right) $\ and its generating curve becomes
either%
\begin{equation*}
\left( 0,u,c_{1}J_{0}\left( \sqrt{\lambda }u\right) +c_{2}Y_{0}\left( \sqrt{%
\lambda }u\right) \right) 
\end{equation*}%
or%
\begin{equation*}
\left( 0,u,c_{1}I_{0}\left( \sqrt{-\lambda}u\right) +c_{2}K_{0}\left( \sqrt{-
\lambda }u\right) \right) ,\text{ }c_{1},c_{2}\in 
\mathbb{R}
.
\end{equation*}
\end{theorem}

\begin{remark} Theorem 4.2 asserts that a surface of revolution in $\mathbb{SI}%
^{3}$\ that satisfy $\bigtriangleup r_{i}=\lambda _{i}r_{i}$ is of null
2-type.
\end{remark}

\begin{example} Consider the surface of revolution $M$ in $\mathbb{SI}^{3}$ given
by%
\begin{equation*}
\left( u\sinh v,u\cosh v,J_{0}\left( u\right) \right) ,\text{ }\left( u,v\right) \in \left[
1,4\right] \times \left[ -1,1\right].
\end{equation*}%
Then $M$ is of null 2-type and thus we draw it and its profile curve as in Fig 2. 
\end{example}

\begin{minipage}{0.55\textwidth}
\includegraphics[scale=0.4]{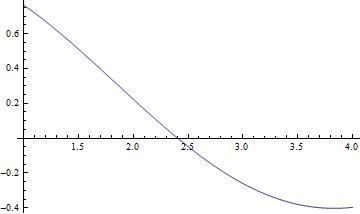} 
\verb Figure2(a).  Profile curve $\left( 0,u,J_{0}\left( u\right) \right)$.
\end{minipage}
\begin{minipage}{0.4\textwidth}
\includegraphics[scale=0.4]{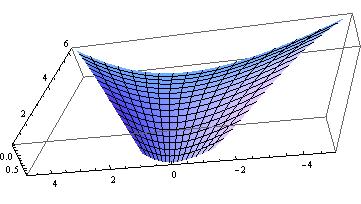}
\verb Figure2b.  A surface of revolution of null 2-type.
\end{minipage}

\bigskip

Next suppose that $M$ has no parabolic points. Then $LN-M^{2}=-f^{\prime
}f^{\prime \prime }u\neq 0,$ i.e. $f$ is a non-linear function. In terms of
the Laplace operator $\bigtriangleup ^{II}$ we get%
\begin{equation}
\left\{ 
\begin{array}{l}
\bigtriangleup ^{II}\left( r_{1}\left( u,v\right) \right) =B\left( u\right)
\sinh v-\frac{1}{f^{\prime }\left( u\right) }\sinh v, \\ 
\bigtriangleup ^{II}\left( r_{2}\left( u,v\right) \right) =B\left( u\right)
\cosh v-\frac{1}{f^{\prime }\left( u\right) }\cosh v, \\ 
\bigtriangleup ^{II}\left( r_{3}\left( u,v\right) \right) =B\left( u\right)
f^{\prime }\left( u\right)+1,%
\end{array}%
\right.  \tag{4.7}
\end{equation}%
where 
\begin{equation*}
B\left( u\right) =\frac{1}{2f^{\prime \prime }\left( u\right) }\left( \frac{%
f^{\prime }\left( u\right) +uf^{\prime \prime }\left( u\right) }{f^{\prime
}\left( u\right) u}-\frac{f^{\prime \prime \prime }\left( u\right) }{%
f^{\prime \prime }\left( u\right) }\right) . \tag{4.8}
\end{equation*}%
Taking $\bigtriangleup ^{II}r_{i}=\lambda _{i}r_{i},$ $\lambda _{i}\in 
\mathbb{R}
,$ $i=1,2,3,$ in $\left( 4.7\right) $ gives%
\begin{equation}
\left\{ 
\begin{array}{l}
B\left( u\right) -\frac{1}{f^{\prime }\left( u\right) }=\lambda u \\ 
B\left( u\right) f^{\prime }\left( u\right) +1=\mu f.%
\end{array}%
\right.  \tag{4.9}
\end{equation}%
where $\lambda _{1}=\lambda _{2}=\lambda $ and $\lambda _{3}=\mu .$ In order
to solve (4.9) we have to distinguish several cases.

\bigskip

\textbf{Case 1. }$\lambda =\mu =0.$ This immediately yields a contradiction
from $\left( 4.9\right) .$

\bigskip 

\textbf{Case 2. }$\lambda =0$ and $\mu \neq 0.$ From first line of $\left(
4.9\right) ,$ we get $B\left( u\right) =\frac{1}{f^{\prime }\left( u\right) 
}$ and considering it into the second line of $\left( 4.9\right) $ we obtain
a contradiction since $f$ must be a non-linear function.

\bigskip

\textbf{Case 3. }$\lambda \neq 0$ and $\mu =0.$ The second line of $\left(
4.9\right) $ gives $B\left( u\right) =\frac{-1}{f^{\prime }\left( u\right) }$
and putting it into the first line of $\left( 4.9\right) $ yields%
\begin{equation*}
f\left( u\right) =\frac{-2}{\lambda }\ln u+c,\text{ }c\in 
\mathbb{R}
.
\end{equation*}%
Note that in this case the surface of revolution $M$ is (s-i)-minimal.

\bigskip

\textbf{Case 4.} $\lambda \neq 0\neq \mu .$ It follows from $\left(
4.9\right) $ that%
\begin{equation}
\lambda uf^{\prime }-\mu f=-2.  \tag{4.10}
\end{equation}%
By solving $\left( 4.10\right) $ we derive%
\begin{equation}
f\left( u\right) =\frac{2}{\mu }+cu^{\frac{\mu }{\lambda }},\text{ }c\in 
\mathbb{R}
.  \tag{4.11}
\end{equation}%
Substituting $\left( 4.11\right) $ into $\left( 4.8\right) $ implies%
\begin{equation}
B\left( u\right) =\frac{u^{\frac{-\mu }{\lambda }+1}}{c\frac{\mu }{\lambda }%
\left( \frac{\mu }{\lambda }-1\right) }.  \tag{4.12}
\end{equation}%
By considering $\left( 4.11\right)$ and $\left( 4.12\right) $ into the second line of $\left(
4.9\right),$ we conclude again a contradiction.

Therefore we have proved the following results.

\begin{theorem} Let $M$\ be a surface of revolution in $\mathbb{SI}^{3}$%
\ given by $\left( 4.1\right) .$\ If $\bigtriangleup ^{II}r_{i}=\lambda
_{i}r_{i}$\ holds then $\left( \lambda _{1},\lambda _{2},\lambda _{3}\right)
=\left( \lambda ,\lambda ,0\right) ,$\ $\lambda \neq 0,$\ and its generating
curve is of the form%
\begin{equation*}
\left( 0,u,\frac{2}{\lambda }\ln u+c\right), \text{ } c \in\mathbb{R} .
\end{equation*}
\end{theorem}

Theorem 4.1 and Theorem 4.3 immediately yields the following result.

\begin{corollary} A surface of revolution satisfies $\bigtriangleup
^{II}r_{i}=\lambda _{i}r_{i}$\ in $\mathbb{SI}^{3}$ if and only if it is
(s-i)-minimal.
\end{corollary}

\begin{example} Given the surface of revolution $M$ in $\mathbb{SI}^{3}$ by%
\begin{equation*}
\left( u\sinh v,u\cosh v,\ln u\right) ,\text{ }\left( u,v\right) \in \left[
0.5,5\right] \times \left[ -0.5,1\right] .
\end{equation*}%
Then $M$ satisfies $\bigtriangleup ^{II}r_{i}=\lambda _{i}r_{i}$\ in $%
\mathbb{SI}^{3}$ for $\left( \lambda _{1},\lambda _{2},\lambda _{3}\right)
=\left( -2,-2,0\right) $ and can be drawn as in Fig3.
\end{example}

\begin{minipage}{0.55\textwidth}
\includegraphics[scale=0.45]{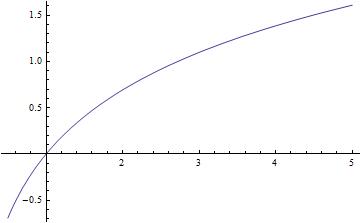} 
\verb Figure3(a).  Profile curve, $f\left(u\right)=\ln u.$ 
\end{minipage}
\begin{minipage}{0.45\textwidth}
\includegraphics[scale=0.35]{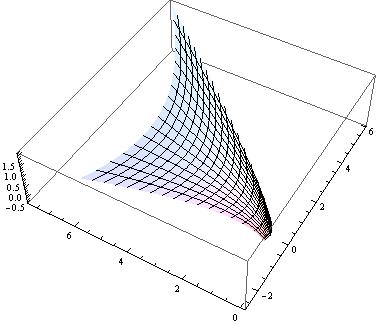}
\verb Figure3b.  A surface of revolution satisfying $\bigtriangleup ^{II}r_{i}=\lambda _{i}r_{i}$.
\end{minipage}

\end{document}